\documentclass[12pt,a4paper]{article}
\usepackage[latin1]{inputenc}
\usepackage[english]{babel}
\usepackage{graphicx}
\usepackage{amsmath}
\usepackage{amsfonts}
\usepackage{amssymb}
\usepackage{indentfirst}
\usepackage{dsfont}
\newtheorem{theorem}{Theorem}[section]
\newtheorem{lemma}[theorem]{Lemma}
\newtheorem{proposition}[theorem]{Proposition}

\newenvironment{definition}[1][Definition]{\begin{trivlist}
\item[\hskip \labelsep {\bfseries #1}]}{\end{trivlist}}

\author{Paul Rochet}
\date{}
\title{A Cram\'er-Rao inequality for non differentiable models}
\begin{document}
\maketitle

\begin{abstract}
We compute a variance lower bound for unbiased estimators in specified statistical models. The construction of the bound is related to the original Cram\'er-Rao bound, although it does not require the differentiability of the model. Moreover, we show our efficiency bound to be always greater than the Cram\'er-Rao bound in smooth models, thus providing a sharper result. 
\end{abstract}

\section{Introduction}

Efficiency theory aims to establish an objective criterion to judge if an estimator is the best possible in a given class. The most famous example is without doubt the Cram\'er-Rao inequality, which states in its simpler form that the variance of an unbiased estimator in a parametric model is not smaller than the inverse of the Fisher information. The inequality was originally stated in \cite{MR0015748} and has been the foundation of a numerous efficiency theories developped in the literature, such as that due to Le Cam and Haj\'ek (see \cite{MR0283911}, \cite{MR0126903}) that extend the Cram\'er-Rao inequality to larger models with alternative regularity assumptions. We refer to \cite{MR1623559} and \cite{MR1652247} for a survey. \bigskip

In this paper, we introduce a variance lower bound for unbiased estimators in a statistical model. The construction of the bound relies on the same idea as the original Cram\'er-Rao bound, although no regularity conditions of any kind are needed. The advantage of our approach is threefold. First, an efficiency bound can be computed without differentiability conditions on the model nor on the parameter to estimate. Second, the bound is adapted to all types of models: parametric, semiparametric or nonparametric. Finally, the efficiency bound is always greater or equal to the Cram\'er-Rao bound (whenever it is well defined) and thus is more informative. \bigskip

The paper is organized as follows. We define our efficiency bound in Section \ref{secconst} and we compare its performance to the Cram\'er-Rao bound in differentiable parametric models. We discuss the generalization to semiparametric models in Section \ref{secsemi} and provide an asymptotic analysis in Section \ref{secasympt}. The proofs of our results are postponed to the Appendix.

\section{Construction of the efficiency bound}\label{secconst}
Let $(\mathcal X,  \mathcal B (\mathcal X)) $ be an open subset of $ \mathbb R^p$ endowed with its borel field, we denote by $\mathcal{P}(\mathcal{X})$ the set of all probability measures on $(\mathcal X,  \mathcal B (\mathcal X)) $. We consider the classical statistical model where we observe an i.i.d. sample $X_1,...X_n$ drawn from an unknown measure $\mu$ and we wish to estimate a parameter $\psi$.\bigskip

The construction of an efficiency bound relies only on two aspects which are the model and the parameter to estimate. The model is defined as the set of possible values for the measure $\mu$. We shall assume in the sequel that the model is well chosen so that $\mu \in \mathcal P$. A parameter $\psi$ is to be understood as a map $\psi: \mathcal P \to \mathrm H$. In this paper, we restrict to finite dimensional parameters, with $\mathrm H$ a subset of $\mathbb R^p$. 

We define the quadratic divergence (or $Q$-divergence) between two probability measures $\mu$ and $\nu$ on $(\mathcal X, \mathcal B)$ as
$$ d(\mu,\nu) = \int_\mathcal X \left( 1 - \textstyle \frac{d\nu}{d\mu} \right)^2 d\mu \  \text{ if } \nu \ll \mu, \ \ d(\mu,\nu) = + \infty \ \text{otherwise}.   $$
The $Q$-divergence is Csisz\'ar's $f$-divergence associated to the convex function $f: x \mapsto (1-x)^2$ (see \cite{MR0219346}). Remark that the $Q$-divergence between two probability measures $\mu$ and $\nu$ is not symmetric, so we shall speak of quadratic divergence of $\nu$ with respect to $\mu$ to avoid confusion. Moreover, let $A$ be a subset of $\mathcal P (\mathcal X )$, we define $d(\mu,A)= \text{inf}_{\nu \in A} \ d(\mu,\nu) $. Any measure $\mu^* \in A$ such that $d(\mu,\mu^*)=d(\mu,A)$ is called \textit{Q-projection} of $\mu$ onto $A$.

\subsection{Main result}\label{secmain}

In the next theorem we show that to each element of a model, can be associated a variance lower bound for an unbiased estimator of a parameter $\psi(\mu) \in \mathbb R^q$. We use the convention $1/\infty = 0$.
 
\begin{theorem}\label{1} Let $\mathcal P$ be a model and $\psi: \mathcal P \to \mathrm H$ a parameter. If $T=T(X_1,...,X_n)$ is an unbiased estimator of $\psi$ in the model $\mathcal P$, then $ \forall \nu \in \mathcal P \setminus \left\lbrace \mu \right\rbrace$:
$$ \emph{var}(T) \geq \ \dfrac{(\psi(\mu) - \psi(\nu))(\psi(\mu) - \psi(\nu))^t}{(d(\mu,\nu) + 1)^n -1}. $$
\end{theorem}
Whenever $\psi$ takes values in $\mathbb R^q$ with $q > 1$, the inequality is meant in the sense of the quadratic forms, i.e. $A \geq B$ if and only if $A - B$ is positive semi-definite. 
Observe that this result does not require any regularity conditions on the model. For instance, it is not needed that $\nu$ be absolutely continuous w.r.t. $\mu$, although in this case the efficiency bound is null and provides no information. \bigskip


\noindent Let $H_\psi^n(\mu,.)$ denote the functional defined on $\mathcal P^* = \mathcal P \setminus \left\lbrace \mu \right\rbrace$ by
$$ H_\psi^n(\mu,\nu) = n \dfrac{(\psi(\mu) - \psi(\nu)) (\psi(\mu) - \psi(\nu))^t}{(d(\mu,\nu) + 1)^n -1}. $$
The quantity $H_\psi^n(\mu,\nu)$ provides a lower bound for $n$ times the variance of an unbiased estimator of $\psi$. Since $H_\psi^n(\mu,\nu)$ is null if $\nu$ is not absolutely continuous w.r.t. $\mu$ or if $ d \nu / d \mu$ is not square $\mu$-integrable, sufficient is to consider the values of $H_\psi^n(\mu,\nu)$ for density measures $\nu = f \mu$ with $f$ in $\mathcal F= \{ f: f \mu \in \mathcal P, \ \int f^2 d \mu < \infty \}$. The main advantage is that $\mathcal F$ being a subspace of $\mathbb L^2(\mu)$, it can be endowed with its natural Hilbert space topology.\bigskip

The result of Theorem \ref{1} gives us all the more information that the right term of the inequality is large. In the case $q >1$, the correct way to interpret this result is to consider real valued linear transformations of $\psi$, where the result can be stated in the form
$$ \forall a \in \mathbb R^q, \ n \ \text{var}(a^t T) \geq  a^t \ H_\psi^n(\mu,\nu)  \ a.  $$
Thus, because the case $q>1$ can be treated by considering real valued parameters, we shall assume for simplicity that $\psi$ takes values in $\mathbb R$, and therefore, $H_\psi^n(\mu,\nu) \in [0; + \infty]$.\bigskip

\noindent We define the efficiency bound for estimating $\psi$ in $\mathcal P$ as the supremum over the whole model
$$ B_\psi^n(\mathcal P) := \underset{\nu \in \mathcal P^*}{\text{sup}} \ H_\psi^n(\mu,\nu)= \underset{f \in \mathcal F \setminus \left\lbrace 1 \right\rbrace}{\text{sup}} \ H_\psi^n(\mu,f \mu) . $$

Let $\Theta $ be a subset of $\mathbb R^d$ and $\{ \mu_\theta \}_{\theta \in \Theta}$ a collection of probability measures on $(\mathcal X,  \mathcal B (\mathcal X)) $ with $\mu_{\theta_0}=\mu$. We say that $\{ \mu_\theta \}_{\theta \in \Theta}$ is \textit{differentiable in }$\mathbb L^2(\mu)$ at $\theta_0$, if there exists a map $g:\mathcal X \rightarrow \mathbb R^d$ such that $\int_\mathcal X g^t g \ d \mu < \infty$ and such that for all $a \in \mathbb R^d$,
$$ \lim_{t\rightarrow 0} \int_\mathcal X \left[ \frac 1 t  \left( \frac{d\mu_{\theta_0 + t a}}{d\mu} (x) -1 \right) - a^t g(x) \right]^2 d \mu(x) = 0.$$ 
The function $g$ is called the score function of the model $\{ \mu_\theta \}_{\theta}$ at $\theta = \theta_0$, while the matrix $\mathcal I = \int_\mathcal X g g^t \ d \mu$ is the Fisher Information. The score can be seen as a Fr\'echet differential of the model $\{ \mu_\theta \}$ in the $\mathbb L^2(\mu)$ sense. More usual definitions of the score generally require the model to be differentiable in an almost-sure sense, which is stronger than the condition above. 
Remark however that, while the differentiability in $\mathbb L^2(\mu)$ is necessary for the sake of this paper, it is less general than the differentiability in quadratic mean, discussed for instance in \cite{MR1915446}.\bigskip

Let $\psi:\{ \mu_\theta \}_{\theta} \to \mathrm H$ be a parameter such that the map $\theta \mapsto \psi(\mu_\theta)$ is differentiable at $\theta_0$ (we note $\dot \psi(\theta_0) \in \mathbb R^{d \times q}$ its derivative matrix), the Cram\'er-Rao inequality states that if $T=T(X_1,...,X_n)$ is an unbiased estimator of $\psi$, then
$$ n \ \text{var} (T) \geq \dot{ \psi}(\theta_0)^t \ \mathcal I^{-1} \ \dot{ \psi}(\theta_0). $$
When the model $\mathcal P$ is differentiable, the Cram\'er-Rao bound $B_\psi(\mathcal P)= \dot{ \psi}(\theta_0)^t \ \mathcal I^{-1} \ \dot{ \psi}(\theta_0)$ provides a variance lower bound for unbiased estimators of $\psi$ in differentiable models. We shall now see in the next proposition the comparison with our efficiency bound $B_\psi^n$ in smooth models.

\begin{proposition}\label{3.2} Let $\{ \mu_\theta \}_{\theta \in \Theta}$ be a differentiable path with $\mu_{\theta_0} = \mu$. Let $\psi: \{ \mu_\theta \}_\theta \rightarrow \mathbb R$ be a map such that $\theta \mapsto \psi(\mu_\theta)$ is differentiable at $\theta_0$
. Then, $B_\psi(\{ \mu_\theta \}_\theta) = \lim_{ \theta \rightarrow \theta_0} H_\psi^n (\mu,\mu_\theta)$ for all $n \in \mathbb N$. In particular,
$$ B_\psi^n( \{ \mu_\theta \}_\theta) \geq B_\psi(\{ \mu_\theta \}_\theta).$$
\end{proposition}
The efficiency bound $B_\psi^n$ improves on the Cram\'er-Rao bound since it is defined as the supremum of $ \nu \mapsto H_\psi^n (\mu,.)$ on the model, while $B_\psi$ is the limit at $\nu \to \mu$. As a result, in differentiable models, the functional $H_\psi^n (\mu,.)$ can be extended by continuity at $\mu$ taking the value $H_\psi^n(\mu,\mu) = B_\psi$. In some situations, the two bounds are identical (i.e. the maximum of $H_\psi^n(\mu,\nu)$ is reached as $\nu \to \mu$), for example as soon as the Cram\'er-Rao bound can be reached for finite samples. On the other hand, it is not rare to have the strict inequality $B_\psi^n( \{ \mu_\theta \}_\theta) > B_\psi(\{ \mu_\theta \}_\theta)$, as we show in the following examples.\\

\textit{Example 1 (Gaussian model).} Consider the Gaussian model $\{ \mu_\theta \}_{\theta \in \mathbb R}$, where $\mu_\theta \sim \mathcal N(\theta, 1)  $ and let $\psi : \mu_\theta \mapsto e^\theta$. We take $\mu \sim \mathcal N(0, 1)$ as the distribution of the observations.
In this model, the Cram\'er-Rao bound is $B_\psi = 1$. On the other hand, we have $ d(\mu,\mu_\theta) = e^{\theta^2} - 1$, yielding 
$$   H_\psi^n(\mu,\mu_\theta) = n \ \frac{(1 - e^{\theta})^2}{e^{n \theta^2} - 1}, $$ 
for $\theta \in (-1; + \infty )$. The supremum is reached for $\theta_n=\frac 1 n$, which gives $B_\psi^n = n (e^{1/n} - 1)$. Thus, we observe a strict inequality $B_\psi^n > B_\psi$ for all $ n \in \mathbb N$. In this case, it is interesting to notice that $B_\psi^n$ is the actual variance of the optimal unbiased estimator of $e^\theta$ in this model. \bigskip

\textit{Example 2 (exponentiel model).} Consider the model $\{ \mu_\theta \}_{\theta >0}$, where $\mu_\theta$ is an exponentiel distribution with parameter $\theta$, i.e. $d \mu_\theta (x) = \theta e^{ - \theta x} \mathds 1 \{x \geq 0 \} dx$. We want to estimate the parameter $\psi: \mu_\theta \mapsto \theta$, the true value of the parameter being $\theta_0=1$. Calculation of the Cram\'er-Rao bound gives $B_\psi=1$. On the other hand, the Q-divergence of $\mu_\theta$ w.r.t. $\mu$ is 
$$ d(\mu,\mu_\theta) =  \frac{\theta^2}{2\theta - 1} - 1, \ \text{for } \ \theta > \frac 1 2 \ \text{ and } \  d(\mu,\mu_\theta) = + \infty \ \text{ otherwise}.$$ 
It follows that
$$ H_\psi^n(\mu,\mu_\theta) = \frac{(\theta - 1)^2 (2 \theta - 1)^n}{ \theta^{2n} - (2 \theta - 1)^n} \ \mathds 1 \{ \theta > 1/2 \}.$$ 

\begin{figure}[ht]
	\centering
		\includegraphics[height=7cm,width=7cm]{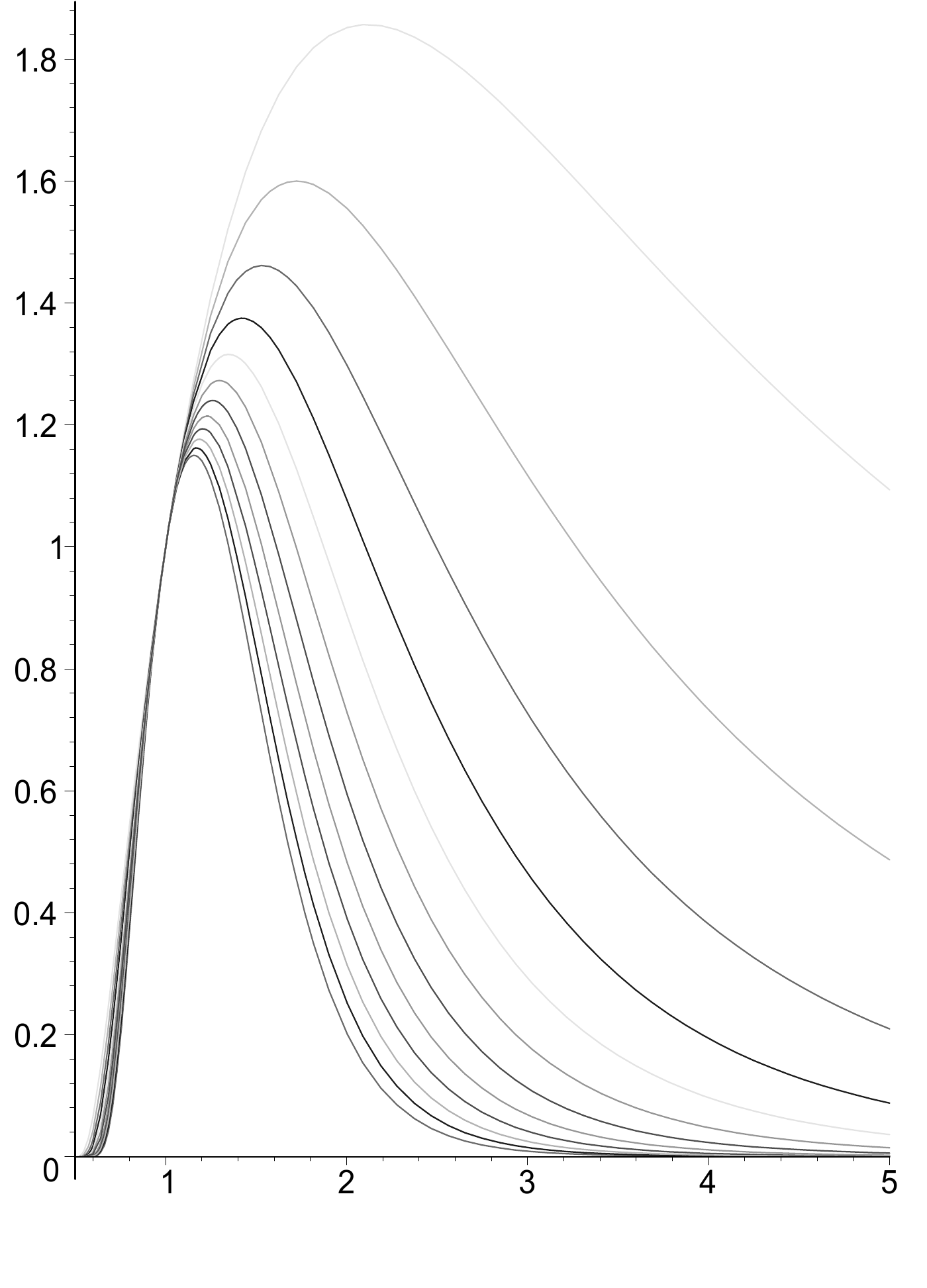}
		\caption{Plot of $\theta \mapsto H_\psi^n(\mu,\mu_\theta)$ for $n=4$ to $15$.}
\end{figure}

The curves are decreasing as $n$ grows (the curve on the top represents $H_\psi^n$ for $n=4$ while the lowest curve is for $n=15$). The functions are not defined at $\theta=1$ but they can be extended by continuity taking the value $H_\psi^n(\mu,\mu) = B_\psi =1$ for all $n \in \mathbb N$. This corresponds on the graph to the intersection point of all the curves. We observe that for all $n \in \mathbb N$, the supremum is larger than the Cram\'er-Rao bound $B_\psi=1$.

\subsection{Application to semiparametric models}\label{secsemi}

Extending the Cram\'er-Rao inequality to semiparametric models can be made using a more general definition of the Fisher Information, calculated by studying differentiable submodels. Based on the idea that, the larger the model, the less information we have, a natural definition of the Fisher Information in large models is to consider the infimum of the Fisher Informations calculated in differentiable submodels (see for instance \cite{}). A \textit{least favorable path} is a differentiable submodel $\{ \mu_\theta \}_{\theta \in \Theta}$ for which the infimum is reached, and therefore, such that $B_\psi(\mathcal P) = B_\psi(\{ \mu_\theta \}_{\theta})$.\bigskip 

The functional $H_\psi^n(\mu,.)$ turns out to be an efficient tool to construct a least favorable path. To see it, consider the level sets $ \mathcal F_\theta = \left\lbrace f \in \mathbb L^2(\mu): \ \psi(f \mu) = \theta  \right\rbrace$ for all values $\theta$ taken by the parameter $\psi$. Setting $\theta_0 = \psi(\mu)$, the expression of the efficiency bound can be written as 
\begin{eqnarray}\label{eq1} B_\psi^n = \underset{\theta \neq \theta_0}{\text{sup}} \ \underset{f \in \mathcal F_\theta }{\text{sup}} \ H_\psi^n(\mu,f \mu) = \sup_{\theta \neq \theta_0} \ n \ \frac{(\theta - \theta_0)^2}{(d(\mu, \mathcal F_\theta)+1)^n-1}.     
\end{eqnarray}
In these settings, we see that calculating the efficiency bound is reduced to maximizing a function of $\theta$. 
The idea is that if we choose the least favorable density in each set $ \mathcal F_\theta$, that is, a function $f_\theta$ maximizing $f \mapsto H_\psi^n(\mu,f \mu)$, the resulting submodel would have to be a least favorable path (if a least favorable measure can not be reached, we may consider a proper collection of densities arbitrarily close to the least favorable measure in each set $ \mathcal F_\theta$, leading to a collection of submodels). Since by construction, the term $\psi(f \mu) - \psi(\mu)$ is constant when $f$ ranges over $ \mathcal F_\theta$, a density maximizing $H_\psi^n(\mu,.)$ on $ \mathcal F_\theta$ is in fact a minimizer of $f \mapsto d(\mu,f \mu)$, which explains the term $d(\mu, \mathcal F_\theta)$ in \eqref{eq1}. 

\begin{definition} We call \textit{quadratic projection path} (or \textit{Q-projection path}) a submodel $\{ \mu_\theta \}_{\theta \in \Theta}$ such that $d(\mu,\mu_\theta) = d(\mu, \mathcal F_\theta)$ and $\psi(\mu_\theta)  = \theta$ for all $\theta \in \Theta$.
\end{definition}

A Q-projection path does not necessarily exist, for instance if the infimum of $d(\mu,.)$ on $\mathcal F_\theta$ is not reachable for some values of $\theta$. However, a Q-projection path does exist as soon as the map $f \mapsto \psi(f \mu)$ is continuous on $\mathcal F$ and if $d(\mu,\mathcal F_\theta)$ is finite for all $\theta \in \Theta$. By making this continuity assumption, we avoid considering trivial cases, the efficiency bound being infinite if $f \mapsto \psi(f \mu)$ is not continuous as $f$ tends to $1$. \bigskip

If the sets $ \mathcal F_\theta$ are convex, a Q-projection path $\{ \mu_\theta \}_{\theta \in \Theta}$ is unique, $\mu_\theta$ being defined as the quadratic projection of $\mu$ on $\mathcal P_\theta = \left\{\nu \in \mathcal P: \psi(\nu) = \theta \right\}$. A Q-projection path does not depend on the number of observations, although it contains a maximizer of $H_\psi^n(\mu,.)$ for all $n \in \mathbb N$. In a certain way, it contains the whole information of the model.\bigskip 

As a straightforward consequence of \eqref{eq1}, a Q-projection path $\{ \mu_\theta \}_{\theta}$ satisfies $ B_\psi^n(\mathcal P) = B_\psi^n(\{ \mu_\theta \}_{\theta})$ for all $n \in \mathbb N$. Moreover, remark that a Q-projection path is a least favorable path if and only if it is differentiable at $\mu = \mu_{\theta_0}$. These remarks are illustrated in the following examples.\bigskip

\textit{Example 3 (moment condition model).} Let $\mathcal P = \left\lbrace \nu \in \mathcal P (\mathcal X): \ \int_\mathcal X \Phi d\nu = 0   \right\rbrace$ for $\Phi:\mathcal X \rightarrow \mathbb R^k$ a known map. We want to estimate $\theta_0 = \int h d\mu \in \mathbb R$ where $h \in \mathbb L^2(\mu)$ is a given function. For all $\theta \in \mathbb R$, $\mathcal F_\theta$ is an affine subspace of $\mathbb L^2(\mu)$ of finite dimension, it is therefore closed and convex. Hence, there exists a unique Q-projection path $\{ \mu_\theta \}_{\theta}$, with densities $f_\theta$ w.r.t. $\mu$. Note $h^\perp$ the part of $h$ orthogonal with $\Phi$ in $\mathbb L^2(\mu)$: $h^\perp = h - ( \int h \Phi d\mu )^t$\begin{large}[\end{large}$\int \Phi \Phi^t d \mu$\begin{large}]\end{large}$^{-1} \Phi$, we have:
\begin{eqnarray*}
 f_\theta  = \arg \min _{f \in \mathcal F_\theta} \ \mathbb E(1-f(X))^2 = 1 - \textstyle (\theta_0 - \theta) V^{-1} (h^\perp - \theta)
\end{eqnarray*}
with $ V =  \text{var}(h^\perp(X))$. Moreover, $d(\mu,\mu_\theta) = \mathbb E(1-f_\theta(X))^2 = (\theta_0 - \theta)^2 V^{-1}$, yielding
$$ B_\psi^n =  \underset{\theta \neq \theta_0 }{\text{sup}} \ \dfrac{n(\theta_0 - \theta)^2}{((\theta_0 - \theta)^2 V^{-1} + 1)^n -1} = V. $$
Note that the model $\{ \mu_\theta \}_{\theta \in \mathbb R}$ is smooth, with Cram\'er-Rao bound $B_\psi = B_\psi^n = V$ for all integer $n$.\\

\textit{Example 4 (empirical likelihood).} Assume that the true measure $\mu$ satisfies the constraint $ \int \Phi_{\theta_0} d\mu = 0$ for some known collection of maps $\left\{ \Phi_\theta: \theta \in \Theta \right\}$ and where $\theta_0$ is the parameter we intend to estimate. The sets 
$\mathcal F_\theta = \left\lbrace f \in \mathbb L^2(\mu): \ \int \Phi_{\theta} f d\mu = 0 \right\rbrace$ are closed and convex. Note $\{ \mu_\theta\}_\theta$ the Q-projection path with densities $f_\theta$ given by
$$ f_\theta = \arg \min_{f \in \mathcal F_\theta} \ \mathbb E \left( 1 - f(X) \right)^2 =  1 - \textstyle \left( \int \Phi_\theta d \mu \right)^t \left[ \text{var}( \Phi_\theta(X)) \right]^{-1} \left( \Phi_\theta - \int \Phi_\theta d \mu \right). $$ 
If we assume that $\theta \mapsto \Phi_\theta$ is differentiable in a neighbourhood of $ \theta_0$, with derivative $\nabla \Phi(.)$, the path $\{ \mu_\theta \}_\theta$ is also differentiable and we have 
$$d(\mu,\mu_\theta) = \textstyle \left( \int \Phi_\theta d \mu \right)^t \left[ \text{var}( \Phi_\theta(X)) \right]^{-1} \left( \int \Phi_\theta d \mu \right),$$
yielding
$$ \textstyle B_\psi^n = \underset{\theta  \neq \theta_0 }{\text{sup}} \ \dfrac{n(\theta_0 - \theta)^2}{(d(\mu,\mu_\theta) + 1)^n -1} \underset{n \rightarrow \infty}{\longrightarrow} \left[  \left( \int \nabla \Phi(\theta_0) d \mu \right)^t \left[ \int \Phi_{\theta_0} \Phi_{\theta_0}^t d \mu \right]^{-1} \left( \int  \nabla \Phi(\theta_0) d \mu \right) \right]^{-1}. $$
We recover the asymptotic efficiency bound of \cite{MR1272085} in this model.

\subsection{Asymptotic properties}\label{secasympt}
We are now interested in the asymptotic analysis of the efficiency bound. 
Writing the first order expansion 
$$(d(\mu,\nu) + 1)^n - 1 = n \ d(\mu,\nu) + \frac{n(n-1)}{2} \ d(\mu,\nu)^2 + ...   $$ 
we see that the sequence $\{ H_\psi^n(\mu,.) \}_{n \in \mathbb N} $ is decreasing and converges pointwise toward $0$ as $n \rightarrow \infty$. So, the non negative sequence $ \{ B_\psi^n \}_{n \in \mathbb N}$ is also decreasing and therefore, it converges (or is infinite). We now aim to prove that, in regular situations, the efficiency bound converges toward the Cram\'er-Rao bound.

\begin{lemma}\label{6.0} Assume that $B_\psi^{n_0} < \infty$ for some $n_0 \in \mathbb N$. Then, for all $\varepsilon > 0$, $H_\psi^n(\mu,.)$ converges uniformly towards $0$ on the set $\left\{ \nu \in \mathcal P: d(\mu,\nu) > \varepsilon \right\}$ as $n \rightarrow \infty$.
\end{lemma}

The condition that $B_\psi^{n_0}$ is finite for some integer $n_0$ is necessary to ensure the existence of an unbiased estimator with finite variance, even asymptotically. However, it may occur that this condition is not fulfilled while the Cram\'er-Rao bound exists and is finite.\bigskip

An interpretation of Lemma \ref{6.0} is that for all element $\nu$ of the model with a non zero distance with $\mu$ (so basically any $\nu \in \mathcal P^*$), the increasing number of observations will eventually end up giving too much information so that the true distribution can not be mistaken with $\nu$. Thus, only the behaviour of the measures of the model in the neighborhood of $\mu$ matters asymptotically. As a result, a measure $\nu \in \mathcal P$ far from $\mu$ will no longer have any influence on the variance of an estimator as soon as the number of observations is large enough. 

\begin{theorem}\label{6} Assume that $B_\psi^{n_0}(\mathcal P) < \infty$ for some $n_0 \in \mathbb N$. If there exists a Q-projection path $\{ \mu_\theta \}_{\theta}$ differentiable at $\mu$, then 
$$ \lim_{n \rightarrow \infty} B_\psi^n(\mathcal P) = B_\psi(\mathcal P).$$ 
\end{theorem}

This result is not surprising as we know that the efficiency bound only depends asymptotically on the behaviour of the model in the neighborhood of $\mu$. Remark that the convergence is pointed out in the examples 1 and 2 in Section \ref{secmain}. We emphasize that, in a parametric model $\{ \mu_\theta \}_\theta$, the efficiency bound has a positive limit $B_\psi^\infty$ in non-trivial cases as soon as the map $\theta \mapsto \sqrt{d(\mu,\mu_\theta)}$ is differentiable at $\theta_0$, while the construction of the Cram\'er-Rao bound requires the much stronger condition of differentiability in $\mathbb L^2(\mu)$. Thus, the efficiency bound $B_\psi^n$ is computable in a larger class of models, while providing at least as good an asymptotic analysis as the Cram\'er-Rao inequality in smooth models.

\section{Appendix}

\textit{Proof of Theorem \ref{1}.} First assume that $\psi(\mu) \in \mathbb R$. If $d(\mu,\nu) = + \infty$, the inequality is trivially verified. If not, first remark that
\begin{eqnarray*} 
\psi(\mu) - \psi(\nu) 
 = \textstyle \mathbb E \left( \left( T - \psi(\mu) \right)  \left( 1 - \frac{d\nu^{\otimes n}}{d\mu^{\otimes n}} \right) \right) 
\end{eqnarray*}
where the expectation is meant under the true distribution of the observations, $\mu^{\otimes n}$. Applying Cauchy-Schwarz inequality, we get
$$ \psi(\mu) - \psi(\nu) \leq \sqrt{\text{var} (T) } \textstyle  \sqrt{d(\mu^{\otimes n},\nu^{\otimes n})}.     $$
It is easy to see that $ d(\mu^{\otimes n},\nu^{\otimes n}) = (d(\mu,\nu) + 1)^n -1$, which yields
$$\text{var} (T) \geq  \dfrac{(\psi(\mu) - \psi(\nu))^2}{(d(\mu,\nu) + 1)^n -1} .  $$
If $\psi(\mu) \in \mathbb R^q$ with $q > 1$, we apply the previous result to the estimator $\alpha ^t T \in \mathbb R$ for some $\alpha \in \mathbb R^q$. We get for all $\nu \neq \mu$:
$$   \text{var} (\alpha ^t T) = \alpha^t \text{var} (T) \alpha \geq  \dfrac{(\alpha^t \psi(\mu) - \alpha^t \psi(\nu))^2}{(d(\mu,\nu) + 1)^n -1} = \alpha ^t \dfrac{ ( \psi(\mu) -\psi(\nu)) ( \psi(\mu) - \psi(\nu))^t}{(d(\mu,\nu) + 1)^n -1} \ \alpha.    $$ 
The inequality holds for all $\alpha \in \mathbb R^q$, which proves the result.\bigskip

\vspace{0.5cm}

\textit{ Proof of Proposition \ref{3.2}.} First remark that if $\{ \mu_\theta \}_{ \theta \in \Theta }$ is differentiable in $\mathbb L^2(\mu)$ at $ \mu = \mu_{\theta_0}$ with score $g$, the limit as $\theta \rightarrow \theta_0$ of $d(\mu,\mu_\theta)/(\theta- \theta_0)^2$ exists and is equal to the Fisher information $\int g^2 d \mu$. In particular, we have for a fixed $n \in \mathbb N$, 
$$(d(\mu,\mu_\theta)+1)^n - 1 = n d(\mu,\mu_\theta) + o(\vert \theta - \theta_0 \vert).$$ 
Hence, 
$$ \lim_{\theta \rightarrow \theta_0} H_\psi^n(\mu,\mu_\theta) = \lim_{\theta \rightarrow \theta_0}  \frac{(\psi(\mu) - \psi(\mu_\theta))^2}{(\theta- \theta_0)^2} \ \frac{(\theta- \theta_0)^2}{d(\mu,\mu_\theta)} = B_\psi( \{ \mu_\theta \}_\theta).  $$

\vspace{0.5cm}

\textit{Proof of Lemma \ref{6.0}.} For all $\nu \neq \mu$, we know that $ H_\psi^{n_0}(\mu,\nu) \leq B_\psi^{n_0}$. The sequence $\{ B_\psi^n \}_{n \in \mathbb N}$ is decreasing as $n \rightarrow \infty$, thus, if $n > n_0$, 
$$ \forall \nu \neq \mu, \ H_\psi^n(\mu,\nu) \leq \frac{n B_\psi^{n_0}}{n_0} \ \frac{(d(\mu,\nu) + 1)^{n_0} - 1}{(d(\mu,\nu) + 1)^{n} - 1}.    $$
Since the function $x \mapsto ((x + 1)^{n_0} - 1)/((x + 1)^{n} - 1) $ is decreasing on the interval $(\varepsilon;+\infty)$ as soon as $n \geq n_o(\varepsilon +1)^{n_0}/((\varepsilon +1)^{n_0} -1)$, we conclude that for large enough values of $n$
$$ \forall \varepsilon > 0, \underset{d(\mu,\nu) > \varepsilon}{\text{sup}} \ H_\psi^n(\mu,\nu) \leq  \frac{n B_\psi^{n_0}}{n_0} \ \frac{(\varepsilon + 1)^{n_0} - 1}{(\varepsilon + 1)^{n} - 1}
.   $$ 
The right term tends to $0$ as $n \rightarrow \infty$ for all $\varepsilon > 0$, which ends the proof.\bigskip

\vspace{0.5cm}

\textit{ Proof of Theorem \ref{6}.} The theorem is true if $B_\psi^\infty = 0$. Now, assume that $B_\psi^\infty > 0$, which warrants that $B_\psi^n(\mathcal P) = B_\psi^n( \{ \mu_\theta \}_\theta)$ for all $n \in \mathbb N$. Let $\{ \mu_n \}_{n \in \mathbb N}$ be a sequence of measures in $\{ \mu_\theta \}_\theta$, suitably chosen so that $\lim_{n \rightarrow \infty} H_\psi^n(\mu,\mu_n) = B_\psi^\infty$. 
We want to prove that $\lim_{n \rightarrow \infty} d(\mu,\mu_n) = 0$. By contradiction, if there exists $\varepsilon > 0$ and an increasing sequence of integers $\{ n_k \}_{k \in \mathbb N}$ such that $\forall k \in \mathbb N, d(\mu,\mu_{n_k}) > \varepsilon$, then:
$$ H_\psi^{n_k}(\mu,\mu_{n_k}) \leq  \underset{d(\mu,\nu) > \varepsilon}{\text{sup}} \ H_\psi^{n_k}(\mu,\nu) \stackrel{k \rightarrow \infty}{\longrightarrow} 0  $$
by Lemma \ref{6.0}, which conflicts with the fact that $\lim_{k \rightarrow \infty} H_\psi^{n_k}(\mu,\mu_{n_k}) = B_\psi^\infty > 0$. So, we conclude that $\lim_{n \rightarrow \infty} d(\mu,\mu_n) = 0$. Since $H_\psi^n(\mu,.)$ is pointwise decreasing as $n \rightarrow \infty$, we get that for all $n \in \mathbb N$,
$$ B_\psi^\infty = \lim_{n \rightarrow \infty} H_\psi^n(\mu,\mu_n) \leq \lim_{\theta \rightarrow \theta_0} H_\psi^n(\mu,\mu_\theta) = B_\psi(\{ \mu_\theta \}_\theta). $$
So, $\{ \mu_\theta \}_\theta$ is a least favorable path of the model and therefore satisfies $B_\psi(\mathcal P) = B_\psi(\{\mu_\theta \}_\theta)$, yielding $B_\psi^\infty(\mathcal P) \leq B_\psi(\mathcal P)$. The reverse inequality being an obvious consequence of Proposition \ref{3.2}, we conclude that $B_\psi^\infty(\mathcal P) = B_\psi(\mathcal P)$.

\bibliography{CRboundbib}

\begin{thebibliography}{BKRW98}

\bibitem[BKRW98]{MR1623559}
P.~J. Bickel, C.~Klaassen, Y.~Ritov, and J.~Wellner.
\newblock {\em Efficient and adaptive estimation for semiparametric models}.
\newblock Springer-Verlag, New York, 1998.
\newblock Reprint of the 1993 original.

\bibitem[Csi67]{MR0219346}
I.~Csisz{\'a}r.
\newblock On topology properties of {$f$}-divergences.
\newblock {\em Studia Sci. Math. Hungar.}, 2:329--339, 1967.

\bibitem[H{\'a}j70]{MR0283911}
J.~H{\'a}jek.
\newblock A characterization of limiting distributions of regular estimates.
\newblock {\em Z. Wahrscheinlichkeitstheorie und Verw. Gebiete}, 14:323--330,
  1969/1970.

\bibitem[LC60]{MR0126903}
L.~Le~Cam.
\newblock Locally asymptotically normal families of distributions. {C}ertain
  approximations to families of distributions and their use in the theory of
  estimation and testing hypotheses.
\newblock {\em Univ. california Publ. Statist.}, 3:37--98, 1960.

\bibitem[QL94]{MR1272085}
J.~Qin and J.~Lawless.
\newblock Empirical likelihood and general estimating equations.
\newblock {\em Ann. Statist.}, 22(1):300--325, 1994.

\bibitem[RR45]{MR0015748}
C.~Radhakrishna~Rao.
\newblock Information and the accuracy attainable in the estimation of
  statistical parameters.
\newblock {\em Bull. Calcutta Math. Soc.}, 37:81--91, 1945.

\bibitem[vdV98]{MR1652247}
A.~van~der Vaart.
\newblock {\em Asymptotic statistics}, volume~3 of {\em Cambridge Series in
  Statistical and Probabilistic Mathematics}.
\newblock Cambridge University Press, Cambridge, 1998.

\bibitem[vdV02]{MR1915446}
A.~van~der Vaart.
\newblock Semiparametric statistics.
\newblock In {\em Lectures on probability theory and statistics
  ({S}aint-{F}lour, 1999)}, volume 1781 of {\em Lecture Notes in Math.}, pages
  331--457. Springer, Berlin, 2002.

\end{thebibliography}
\bibliographystyle{alpha}
 
\end{document}